%
\catcode`@=11
%
%
\def\bibn@me{R\'ef\'erences}
\def\bibliographym@rk{\centerline{{\sc\bibn@me}}
	\sectionmark\section{\ignorespaces}{\unskip\bibn@me}
	\bigbreak\bgroup
	\ifx\ninepoint\undefined\relax\else\ninepoint\fi}
%
%
%
\let\refsp@ce=\ 
\let\bibleftm@rk=[
\let\bibrightm@rk=]
%
%
%
\def\numero{n\raise.82ex\hbox{$\fam0\scriptscriptstyle o$}~\ignorespaces}
%
%
\newcount\equationc@unt
\newcount\bibc@unt
\newif\ifref@changes\ref@changesfalse
\newif\ifpageref@changes\ref@changesfalse
\newif\ifbib@changes\bib@changesfalse
\newif\ifref@undefined\ref@undefinedfalse
\newif\ifpageref@undefined\ref@undefinedfalse
\newif\ifbib@undefined\bib@undefinedfalse
\newwrite\@auxout
%
%
\def\eqnum{\global\advance\equationc@unt by 1%
\edef\lastref{\number\equationc@unt}%
\eqno{(\lastref)}}
%
%
%
%
%
%
\def\re@dreferences#1#2{{%
	\re@dreferenceslist{#1}#2,\undefined\@@}}
\def\re@dreferenceslist#1#2,#3\@@{\def\next{#2}%
	\expandafter\ifx\csname#1@@\meaning\next\endcsname\relax
	??\immediate\write16
	{Warning, #1-reference "\next" on page \the\pageno\space
	is undefined.}%
	\global\csname#1@undefinedtrue\endcsname
	\else\csname#1@@\meaning\next\endcsname\fi
	\ifx#3\undefined\relax
	\else,\refsp@ce\re@dreferenceslist{#1}#3\@@\fi}
%
%
%
\def\newlabel#1#2{{\def\next{#1}\newl@bel#2}}
\def\newl@bel#1#2{%
	\expandafter\xdef\csname ref@@\meaning\next\endcsname{#1}%
	\expandafter\xdef\csname pageref@@\meaning\next\endcsname{#2}}
\def\label#1{{%
	\toks0={#1}\message{ref(\lastref) \the\toks0,}%
	\ignorespaces\immediate\write\@auxout%
	{\noexpand\newlabel{\the\toks0}{{\lastref}{\the\pageno}}}%
	\def\next{#1}%
	\expandafter\ifx\csname ref@@\meaning\next\endcsname\lastref%
	\else\global\ref@changestrue\fi%
	\newlabel{#1}{{\lastref}{\the\pageno}}}}
\def\ref#1{\re@dreferences{ref}{#1}}
\def\pageref#1{\re@dreferences{pageref}{#1}}
%
%
\def\bibcite#1#2{{\def\next{#1}%
	\expandafter\xdef\csname bib@@\meaning\next\endcsname{#2}}}
\def\cite#1{\bibleftm@rk\re@dreferences{bib}{#1}\bibrightm@rk}
%
%
\def\beginthebibliography#1{\bibliographym@rk
	\setbox0\hbox{\bibleftm@rk#1\bibrightm@rk\enspace}
	\parindent=\wd0
	\global\bibc@unt=0
	\def\bibitem##1{\global\advance\bibc@unt by 1
		\edef\lastref{\number\bibc@unt}
		{\toks0={##1}
		\message{bib[\lastref] \the\toks0,}%
		\immediate\write\@auxout
		{\noexpand\bibcite{\the\toks0}{\lastref}}}
		\def\next{##1}%
		\expandafter\ifx
		\csname bib@@\meaning\next\endcsname\lastref
		\else\global\bib@changestrue\fi%
		\bibcite{##1}{\lastref}
		\medbreak
		\item{\hfill\bibleftm@rk\lastref\bibrightm@rk}%
		}
	}
\def\endthebibliography{\egroup\par}
%
%
%
\def\@closeaux{\closeout\@auxout
	\ifref@changes\immediate\write16
	{Warning, changes in references.}\fi
	\ifpageref@changes\immediate\write16
	{Warning, changes in page references.}\fi
	\ifbib@changes\immediate\write16
	{Warning, changes in bibliography.}\fi
	\ifref@undefined\immediate\write16
	{Warning, references undefined.}\fi
	\ifpageref@undefined\immediate\write16
	{Warning, page references undefined.}\fi
	\ifbib@undefined\immediate\write16
	{Warning, citations undefined.}\fi}
%
%
\immediate\openin\@auxout=\jobname.aux
\ifeof\@auxout \immediate\write16
  {Creating file \jobname.aux}
\immediate\closein\@auxout
\immediate\openout\@auxout=\jobname.aux
\immediate\write\@auxout {\relax}%
\immediate\closeout\@auxout
\else\immediate\closein\@auxout\fi
%
%
\input\jobname.aux
\immediate\openout\@auxout=\jobname.aux
%
%

\def\bibn@me{R\'ef\'erences bibliographiques}
%
\def\bibliographym@rk{\bgroup}
%
%
\outer\def\bye{ 	\par\vfill\supereject\end}

\def\Z{{\bf Z}} 
\def\R{{\bf R}}

\overfullrule=0pt

\magnification=1200

\def\house#1{\setbox1=\hbox{$\,#1\,$}%
\dimen1=\ht1 \advance\dimen1 by 2pt \dimen2=\dp1 \advance\dimen2 by 2pt
\setbox1=\hbox{\vrule height\dimen1 depth\dimen2\box1\vrule}%
\setbox1=\vbox{\hrule\box1}%
\advance\dimen1 by .4pt \ht1=\dimen1
\advance\dimen2 by .4pt \dp1=\dimen2 \box1\relax}

  \def\eps{{\varepsilon}}

  \def\uy{{\underline y}}

\def\build#1_#2^#3{\mathrel{\mathop{\kern 0pt#1}\limits_{#2}^{#3}}}

\def\date {le\ {\the\day}\ \ifcase\month\or 
janvier\or fevrier\or mars\or avril\or mai\or juin\or juillet\or
ao\^ut\or septembre\or octobre\or novembre\or 
d\'ecembre\fi\ {\oldstyle\the\year}}

\font\fivegoth=eufm5 \font\sevengoth=eufm7 \font\tengoth=eufm10

\newfam\gothfam \scriptscriptfont\gothfam=\fivegoth
\textfont\gothfam=\tengoth \scriptfont\gothfam=\sevengoth

\def\cqfd{\unskip\kern 6pt\penalty 500 \raise 0pt\hbox{\vrule\vbox 
to6pt{\hrule width 6pt \vfill\hrule}\vrule}\par}

\def\smallsquare{\vbox{\hrule\hbox{\vrule height 1 ex\kern 1 ex\vrule}\hrule}}
\def\cqfd{\hfill \smallsquare\vskip 3mm}

\def\hw{{\hat w}} \def\hv{{\hat v}}

\def\cV{{\cal V}}

\def\DD{D_{\theta, \hv}(\beta_0, \beta_1)}

\vskip 2mm

\centerline{\bf Uniform Diophantine approximation related to $b$-ary and $\beta$-expansions}

\vskip 8mm

\centerline{Yann B{\sevenrm UGEAUD} and Lingmin L{\sevenrm IAO} \footnote{}{\rm
}}

{\narrower\narrower
\vskip 12mm

\proclaim Abstract. { 
Let $b\geq 2$ be an integer and $\hat{v}$ a real number.
Among other results, we compute the Hausdorff dimension of the set of real numbers
$\xi$ with the property that, for every sufficiently large
integer $N$, there exists an integer $n$ such that $1 \le n \le N$ 
and the distance between $b^n \xi$ and its nearest integer is at most equal to  
$b^{-\hat{v} N}$. We further solve the same question 
when replacing $b^n\xi$ by $T^n_\beta \xi$, 
where $T_\beta$ denotes the classical $\beta$-transformation.
}

}

\vskip 6mm

\vskip 5mm

\centerline{\bf 1. Introduction and results}

\vskip 6mm

Throughout this text, $\| \cdot \|$ stands for the distance to the
nearest integer. 
Let $\xi$ be an irrational real number.
The well-known Dirichlet Theorem asserts that for 
{\it every} real number $X \ge 1$, there exists
an integer $x$ with $1 \le x \le X$ and
$$
\| x  \xi \|  < X^{-1}.  \eqno (1.1)
$$
This is a uniform statement in the sense that (1.1) has a solution for
every sufficiently large real number $X$ (as opposed to `for arbitrarily large real
numbers $X$').
Following the notation introduced in \cite{BuLa05a}, we denote by $\hw_1(\xi)$
the supremum of the real numbers $w$ such that, for
any sufficiently large real number $X$, the inequality
$$
\| x  \xi \| < X^{-w} 
$$
has an integer solution $x$ with $1 \le x \le X$. 
The Dirichlet Theorem implies that $\hw_1 (\xi) \ge 1$. 
In 1926, Khintchine \cite{Kh26} 
established that, in fact, $\hw_1 (\xi) = 1$ always holds. 
To see this, let $(p_\ell / q_\ell)_{\ell \ge 1}$
denote the sequence of convergents to $\xi$. If $\hw_1 (\xi) > 1$, then there exists a
positive real number $\eps$ such that, for every sufficiently large integer $\ell$, the
inequality $\| q \xi \| < (q_\ell-1)^{-1 - \eps}$ has an integer solution  
$q$ with $1 \le q < q_\ell$.
However, it follows from the theory of continued fractions that $\| q \xi \| \ge \|q_{\ell-1} \xi\|
\ge 1/ (2 q_\ell)$. This gives $q_{\ell}^{\eps/2} < 2$, thus we have reached a contradiction. 
Consequently, the set of values taken by the exponent  
of Diophantine approximation $\hw_1$ is easy to determine.    

In the present paper, we first restrict our attention to approximation by rational 
numbers whose denominator is a power of some given integer $b \ge 2$
and we consider the following exponents of approximation.

\proclaim Definition 1.1.
Let $\xi$ be an irrational real number.
Let $b$ be an integer with $b \ge 2$.
We denote by $v_b (\xi)$ the supremum
of the real numbers $v$ for which the equation
$$
\| b^n \xi \| < (b^n)^{-v}
$$
has infinitely many solutions in positive integers $n$.
We denote by $\hv_b (\xi)$ the supremum
of the real numbers $\hv$ for which, for every sufficiently large
integer $N$, the equation
$$
\| b^n \xi \|< (b^N)^{-\hv}
$$
has a solution $n$ with $1 \le n \le N$.

The exponents $v_b$ have already been introduced in \cite{AmBu10}; see also 
Chapter 7 of \cite{BuLiv2}.
Roughly speaking, the quantity $v_b(\xi)$ measures the maximal lengths of blocks  
of digits $0$ (or of digits $b-1$) in the $b$-ary expansion of $\xi$. 
The exponents $\hv_b$ are, like $\hw_1$, exponents of uniform   
approximation. Although they occur rather naturally, they do not seem to   
have been studied until now.   

Alternatively, we may consider the quantity $\ell_n (\xi)$, defined as the maximal 
length of a block of digits $0$ or a block of digits $b-1$ among the $n$ 
first $b$-ary digits of $\xi$. 
Let $v$ and $\hv$ be positive real numbers.
We have $v_b(\xi) \ge v$ (resp., $\hv_b(\xi) \ge \hv$) if, and only if, 
there are arbitrarily large integers $n$ (resp., for every sufficiently large
integer $n$) such that $\ell_n (\xi)/n \ge v / (1 + v)$
(resp., $\ell_n (\xi)/n \ge \hv / (1 + \hv)$).

It follows immediately from the above mentioned result of Khintchine that
$\hv_b (\xi) \le 1$. 
An easy covering argument shows that the set
$$
\{\xi \in \R : v_b (\xi) = 0 \} 
$$
has full Lebesgue measure. Since
$$
0 \le\hv_b (\xi)\le  v_b (\xi), \eqno (1.2)
$$
almost all real numbers $\xi$ (with respect to the Lebesgue measure)
satisfy $\hv_b (\xi) = 0$.
However, it is easy to construct a suitable lacunary series $f(x)$   
such that $\hv_b (f(1/b))$ has a prescribed value between $0$ and $1$.  
Indeed, for any $v > 0$ we have
$$
v_b \Bigl( \sum_{j \ge 1} b^{-(1+v)^j} \Bigr) = v
\quad
\hbox{and}
\quad
\hv_b \Bigl( \sum_{j \ge 1} b^{-(1+v)^j} \Bigr) =  {v \over v + 1}.
$$
Observe also that
$$
v_b \Bigl( \sum_{j \ge 1} b^{-2^{j^2}} \Bigr) = + \infty
\quad
\hbox{and}
\quad
\hv_b \Bigl( \sum_{j \ge 1} b^{-2^{j^2}} \Bigr) = 1.
$$
Hence, for every $v$ in $\R_{>0} \cup \{+ \infty\}$
and every real number $\hv$ in $[0, 1]$, the sets
$$
\cV_b (v) := \{\xi \in \R : v_b (\xi) = v \} 
\quad
\hbox{and}
\quad
\widehat{\cV_b} (\hv) := \{\xi \in \R : \hv_b (\xi) = \hv \} 
$$
are non-empty. Again, an easy covering argument yields that
$$
\dim \cV_b (+ \infty) = 0,  \eqno (1.3) 
$$
where $\dim$ stands for the Hausdorff dimension.

Let $\hv$ be in $[0, 1]$. It follows from (1.2) that  
$$
\dim \widehat{\cV_b}  (\hv) 
\le \dim \{\xi \in \R : \hv_b (\xi) \ge  \hv \} 
\le \dim \{\xi \in \R : v_b (\xi) \ge  \hv \}.
$$
Combined with
$$
\dim \{\xi \in \R : v_b (\xi) \ge  \hv \} = {1 \over 1 + \hv}, \eqno (1.4) 
$$
which follows from a general result of Borosh and Frankel \cite{BoFr72},
this gives
$$
\dim \widehat{\cV_b}  (\hv)  \le  {1 \over 1 + \hv}.   
$$
Note that (1.4) is also a special case of Theorem 5 in \cite{PeSc08} and that, 
furthermore, it easily  
follows from the mass transference principle of Beresnevich and Velani \cite{BeVe06}. 
Moreover, by \cite{BeVe06} or by Theorem 7 from \cite{AmBu10}, we have 
$$
\dim \cV_b (v) = {1 \over 1 + v}, \eqno (1.5)
$$
for every $v \ge 0$.  

Our first result gives the Hausdorff dimension of the set $\widehat{\cV_b} (\hv)$ for 
any $\hv$ in $[0, 1]$. 

\proclaim Theorem 1.1.
Let $b \ge 2$ be an integer and $\hv$ be a real number in $[0, 1]$.
Then we have
$$
\dim \{\xi \in \R : \hv_b (\xi) \ge  \hv \} =
\dim \{\xi \in \R : \hv_b (\xi) = \hv \} =
\Bigl( {1-\hv \over 1+\hv} \Bigr)^2.    \eqno (1.6)
$$

Theorem 1.1 follows from a more general statement, in which the
values of both functions $v_b$ and $\hv_b$ are prescribed.

\proclaim Theorem 1.2.
Let $b \ge 2$ be an integer. 
Let $\theta$ and $\hv$ be positive real numbers with $\hv < 1$.
If $\theta < 1/(1 - \hv)$, then the set
$$
\{ \xi \in \R : \hv_b (\xi) \ge \hv\} \cap \{ \xi \in \R : v_b (\xi) = \theta \hv \}
$$
is empty. Otherwise, we have
$$
\dim (\{ \xi \in \R : \hv_b (\xi) = \hv\} \cap \{ \xi \in \R : v_b (\xi) = \theta \hv \}) =
{\theta - 1 -  \theta \hv  \over (1 + \theta \hv) ( \theta - 1) }.    \eqno (1.7)
$$
Furthermore,
$$
\dim \{ \xi \in \R : \hv_b (\xi) =  1\} = 0. \eqno (1.8)
$$

A key observation in the proof of Theorem 1.2 is the fact that
the right hand side inequality of (1.2) can be considerably improved.
Namely, we show in Subsection 2.1 that
$$
v_b (\xi)  \enspace \hbox{is infinite when $\hv_b(\xi) = 1$}   \eqno (1.9) 
$$
and
$$
v_b(\xi) \ge \hv_b(\xi) / (1 - \hv_b(\xi)) \enspace \hbox{when $\hv_b(\xi) < 1$}. 
$$
The latter inequality immediately implies the first statement of Theorem 1.2.  
Furthermore, the combination of (1.3) and (1.9) gives (1.8). 

A rapid calculation shows that
the right hand side of (1.7) is a continuous function of the parameter $\theta$ 
on the interval $[1/(1 - \hv), + \infty)$,  
reaching its maximum at the point $\theta_0 := 2/(1 - \hv)$ and only at that point.  
This maximum is precisely equal to $(1 - \hv)^2 / (1 + \hv)^2$, 
namely the right hand side of (1.6). 
This essentially shows that Theorem 1.2 implies Theorem 1.1    
(a complete argument is given at the end of Subsection 2.1).    

We remark that Theorem 1.2 allows us to reprove (1.5). 
To see this, write $v=\theta \hv$, then $\hv=v/\theta$ and (1.7) becomes
$$
\eqalign{
\dim (\{ \xi \in \R : \hv_b (\xi) = v/\theta \} \cap \{ \xi \in \R : v_b (\xi) =  v \}) 
& = {\theta - 1 -  v  \over (1 + v) ( \theta - 1) } \cr
& = {1 \over 1+v} \Bigl(1-{v \over \theta-1} \Bigr).  \cr}
$$
Letting $\theta$ tend to infinity, we see that $\dim \cV_b(v) \ge 1 / (1 + v)$.  
The reverse inequality can easily be obtained by using the natural covering.    

\bigskip
Beside $b$-ary expansions, we can as well consider $\beta$-expansions.
For $\beta > 1$, let $T_{\beta}$ be the $\beta$-transformation defined on $[0, 1]$ by 
$$ 
T_\beta(x):=\beta x-\lfloor \beta x \rfloor, 
$$
where $\lfloor \cdot \rfloor$ denotes the integer part function.
We assume that the reader is familiar with the classical results on $\beta$-expansions. 
Some useful facts are recalled in Section 3.

For a real number $\beta > 1$, 
we define in a similar way the functions $v_{\beta}$ and $\hv_{\beta}$.

\proclaim Definition 1.3.
Let $\beta > 1$ be a real number. 
Let $x\in [0,1]$.
We denote by $v_\beta (x)$ the supremum
of the real numbers $v$ for which the equation
$$
T_\beta^n x  < (\beta^n)^{-v}
$$
has infinitely many solutions in positive integers $n$.
We denote by $\hv_\beta (x)$ the supremum
of the real numbers $\hv$ for which, for every sufficiently large
integer $N$, the equation
$$
T_\beta^n x  < (\beta^N)^{-\hv}  
$$
has a solution $n$ with $1 \le n \le N$.

Observe that Definitions 1.1 and 1.3 do not coincide when $\beta$ is an integer
at least equal to $2$. However, this should not cause any trouble and we may consider
that the following result, established by
Shen and Wang \cite{ShWa13}, extends (1.4) to $\beta$-expansions. 

\proclaim Theorem SW. 
Let $\beta > 1$ be a real number and $v$ be a positive real number. Then,
$$
\dim \{x \in [0,1] : v_{\beta} (x) \ge v\} = {1 \over 1 + v}.
$$

We establish the following analogues of Theorems 1.1 and 1.2 
for $\beta$-expansions.

\proclaim Theorem 1.4.
Let $ \beta > 1$ be a real number. 
Let $\theta$ and $\hv$ be positive real numbers with $\hv < 1$.
If $\theta < 1/(1 - \hv)$, then the set
$$
\{ x \in [0,1]  : \hv_\beta (x) \ge \hv\} \cap \{ x \in [0,1] : v_\beta (x) = \theta \hv \}
$$
is empty. Otherwise, we have
$$
\dim (\{ x \in [0,1] : \hv_\beta (x) = \hv\} \cap \{ x \in [0,1] : v_\beta (x) = \theta \hv \}) =
{\theta - 1 -  \theta \hv  \over (1 + \theta \hv) ( \theta - 1) }.       
$$
Furthermore,
$$
\dim (\{ x \in [0,1] : \hv_\beta (x) =  1\}) = 0.
$$

In the same way as Theorem 1.2 implies Theorem 1.1, the next statement follows from 
Theorem 1.4.

\proclaim Theorem 1.5.
Let $\beta > 1$ be a real number and $\hv$ be a real number in $[0, 1]$.
Then we have
$$
\dim \{x \in [0,1] : \hv_{\beta} (x) \ge \hv \} = 
\dim \{x \in [0,1] : \hv_{\beta} (x) = \hv \} =
\Bigl( {1- \hv \over 1+ \hv} \Bigr)^2.
$$

Persson and Schmeling \cite{PeSc08} have adopted another point of view, by considering
the $\beta$-expansions of $1$ and letting $\beta$ vary.

\proclaim Theorem PS.
Let $\beta_0, \beta_1$ and $v$ be real numbers with $1 < \beta_0 < \beta_1 < 2$
and $v > 0$. Then,
$$
\dim \{\beta \in (\beta_0 , \beta_1) : v_{\beta} (1) \ge v \} = {1 \over 1 + v}.
$$

The assumption $\beta_1 < 2$ in Theorem PS can easily be removed; see \cite{LPWW}. 
Applying some of the ideas from \cite{PeSc08}, 
we obtain the following theorem. 

\proclaim Theorem 1.6. 
Let $\theta$ and $\hv$ be positive real numbers with $\hv < 1$.
If $\theta < 1/(1 - \hv)$, then the set
$$
\{ \beta>1 : \hv_\beta (1) \ge \hv\} \cap \{ \beta>1 : v_\beta (1) = \theta \hv \}
$$
is empty. Otherwise, we have
$$
\dim (\{ \beta > 1 : \hv_\beta (1) = \hv\} \cap \{ \beta > 1 : v_\beta (1) = \theta \hv \}) =
{\theta - 1 -  \theta \hv  \over (1 + \theta \hv) ( \theta - 1) }.    
$$
Furthermore,
$$
\dim \{\beta > 1 : \hv_{\beta} (1) \ge \hv \} 
= \Bigl( {1- \hv \over 1+ \hv} \Bigr)^2, 
$$
and
$$
\dim (\{ \beta > 1 : \hv_\beta (1) =  1\}) = 0. 
$$


Our paper is organized as follows. Theorem 1.2 is proved in Section 2.
We recall classical results from the theory of $\beta$-expansion in Section 3 and establish
Theorems 1.4 and 1.6 in Sections 4 and 5, respectively. 
Diophantine approximation on Cantor sets is briefly discussed in Section 6.


Throughout this text, we denote by $|I|$ the length of the interval $I$
and we use $\sharp$ to denote the cardinality of a finite set.

\vskip 5mm

\centerline{\bf 2. Proof of Theorem 1.2}

\vskip 6mm

\noindent {\bf 2.1. Upper bound.}

Let $\hv$ be a real number with $0 < \hv \le 1$.
We wish to bound from above the dimension of
$$
\{ \xi \in \R : \hv_b (\xi) \ge \hv\}.
$$
By (1.3), it is sufficient to consider the set   
$$
\{ \xi \in \R : \hv_b (\xi) \ge \hv\} \cap \{ \xi \in \R : v_b (\xi) < + \infty \}.
$$
Let $\xi$ be an irrational real number. 
Throughout this section, in view of the preceding observations,
we assume that 
$$
0 < \hv_b(\xi) \le  v_b (\xi) < + \infty.
$$ 

Let
$$
\xi := \lfloor \xi \rfloor + \sum_{j \ge 1} {a _j \over b^j}
$$
denote the $b$-ary expansion of $\xi$. It is understood that the digits $a_1, a_2, \ldots$
all belong to the set $\{0, 1, \ldots , b-1\}$.

Define the increasing sequences $(n'_k)_{k \ge 1}$ and $(m'_k)_{k \ge 1}$
as follows: for $k \ge 1$, we have either
$$
a_{n'_k} >0, a_{n'_k + 1} = \ldots = a_{m'_k - 1} = 0, a_{m'_k} > 0
$$
or
$$
a_{n'_k} < b-1, a_{n'_k + 1} = \ldots = a_{m'_k - 1} = b-1, a_{m'_k} < b-1.
$$
Furthermore, for every $j$ such that $a_j = 0$ or $b-1$, there exists an index $k$
satisfying $n'_k < j < m'_k$. Since $v_b (\xi)$ is positive, we get  
$$
\limsup_{k \to + \infty} \, (m'_k - n'_k) = + \infty.  \eqno (2.1) 
$$
Now, we take the maximal subsequences $(n_k)_{k \ge 1}$ and $(m_k)_{k \ge 1}$
of $(n'_k)_{k \ge 1}$ and $(m'_k)_{k \ge 1}$, respectively, in such a way that the sequence 
$(m_k - n_k)_{k \ge 1}$ is non-decreasing. More precisely, take $n_1=n'_1$ and $m_1=m'_1$. 
Let $k \ge 1$ be such that $n_k=n'_{j_k}$ and $m_k=m'_{j_k}$ have been defined. Set  
$$
j_{k+1}:=\min \{j>j_k: m'_j-n'_j \geq m_k-n_k\}.
$$
Then, define
$$
n_{k+1}=n'_{j_{k+1}}\quad {\rm and} \quad m_{k+1}=m'_{j_{k+1}}.
$$
Observe that, by (2.1), the sequence $(j_k)_{k \ge 1}$ is well defined. 
Furthermore, $m_k-n_k$ tends to infinity as $k$ tends to infinity. 

Note that
$$
b^{n_k - m_k} < \| b^{n_k} \xi \| < b^{n_k - m_k + 1}.
$$
By construction, we have
$$
v_b (\xi) = \limsup_{k \to + \infty} \, {m_k - n_k \over n_k} =   
\limsup_{k \to + \infty} \, {m_k  \over n_k} - 1     \eqno (2.2)
$$
and
$$
\hv_b (\xi) = \liminf_{k \to + \infty} \, {m_k - n_k  \over n_{k+1}}   
\le \liminf_{k \to + \infty} \, {m_k - n_k \over m_k} = 
1 - \limsup_{k \to + \infty} \, {n_k  \over m_k}.   \eqno (2.3)
$$
Since
$$
\Bigl(\limsup_{k \to + \infty} \, {n_k  \over m_k}\Bigr) \cdot 
\Bigl(\limsup_{k \to + \infty} \, {m_k  \over n_k}\Bigr) \ge 1,
$$
we derive from (2.2) and (2.3) that   
$$
\hv_b (\xi) \le 1 - {1 \over 1 + v_b (\xi)} = {v_b (\xi) \over 1 + v_b (\xi)}.
$$
Noticing that $\hv_b (\xi) < 1$ since $v_b (\xi)$ is assumed to be finite, 
we have proved that
$$
v_b (\xi) \ge {\hv_b (\xi) \over 1 - \hv_b (\xi)}.   \eqno (2.4)
$$

Let $\xi$ be a real number with $\hv_b (\xi) \ge \hv$.  
Take a subsequence $(k_j)_{j \ge 1}$ along which the supremum of (2.2) is obtained. 
For simplicity, we still write $(n_k)_{k \ge 1}, (m_k)_{k \ge 1}$  
for the subsequences $(n_{k_j})_{j \ge 1}, (m_{k_j})_{j \ge 1}$.   
We remark that when passing to the subsequence, 
the first equality in (2.3) becomes an inequality.   

Let $\eps$ be a real number with $0 < \eps < v_b(\xi) / 2$.  
Observe that for $k$ large enough we have
$$
(v_b(\xi) - \eps) n_k\le m_k-n_k  \le (v_b(\xi) + \eps) n_k    \eqno (2.5)   
$$
and
$$
m_k -n_k \ge n_{k+1} (\hv - \eps).  \eqno (2.6)   
$$
The last inequality means that the length of the block of $0$ (or $b-1$)
starting at index $n_k+1$ is at least equal to $n_{k+1} (\hv - \eps)$. 

The combination of the second inequality of (2.5) and (2.6) gives
$$
(v_b(\xi) + \eps) n_k \ge (\hv - \eps) n_{k+1}.
$$
Consequently, there exist an integer $n'$ and a positive real number $\eps'$ such that 
the sum of all the lengths of the blocks of $0$ or $b-1$
in the prefix of length $n_k$ of the infinite sequence $a_1 a_2 \ldots$ is,
for $k$ large enough, at least equal to
$$
\eqalign{
(\hv - \eps) n_k \Bigl( 1 + {\hv - \eps \over v_b(\xi) + \eps} +
\Bigl( {\hv - \eps \over v_b(\xi) + \eps} \Bigr)^2 + \ldots \Bigr) - n' 
& =  n_k \, {(\hv - \eps) (v_b(\xi) + \eps) \over v_b(\xi) - \hv + 2 \eps} - n' \cr
& \ge n_k \, \Bigl( {\hv  v_b(\xi)  \over v_b(\xi) - \hv} - \eps' \Bigr) \cr
& = n_k \, \Bigl( {\theta \hv \over \theta - 1} - \eps' \Bigr),  \cr} 
\eqno (2.7)
$$
where the parameter $\theta$ is defined by
$$
v_b (\xi) = \theta \hv.
$$
Note that, by (2.4), we must have
$$
\theta \ge {1 \over 1 - \hv}.  \eqno (2.8)
$$

By the first inequality of (2.5), we have 
$$
m_k\ge (1+v_b(\xi)-\eps)n_k \ge (1+v_b(\xi)-\eps)m_{k-1},
$$
for $k$ large enough. 
Thus $(m_k)_{k \ge 1}$ increases at least exponentially.  
Since $n_k\geq m_{k-1}$ for $k \ge 2$, the sequence $(n_k)_{ k \ge 1}$  
also increases at least exponentially. 
Consequently, there exists a positive real number $C$ such that $k\le C \log n_k$, 
for $k$ large enough. 

Now, let us construct a covering. 
Remind that all the integers $n_k, m_k$ defined above depend on $\xi$.    
Fix $\hv$ and $\theta$ with $\theta \ge {1 \over 1 - \hv}$. 
Let $(n_k)_{k \ge 1}$ and $(m_k)_{k \ge 1}$ be sequences such that
$$
\lim_{k\to \infty} {m_k-n_k \over n_k} = \theta \hv, \quad {\rm and } 
\quad \liminf_{k\to \infty} {m_k-n_k \over n_{k+1}}\geq \hv.
$$

For fixed $k$, we collect all those $\xi$ with $\hv_b(\xi)=\hv$ 
and $v_b(\xi)=\theta \hv$ whose $b$-ary expansions have 
blocks of $0$ (or of $b-1$) between $n_k$ and $m_k$.

By the precedent analysis, there are at most $C \log n_k$ blocks of $0$ (or of $b-1$)    
in the prefix of length $n_k$ of the infinite sequence $a_1 a_2 \ldots$.
Since there are, obviously, at most $n_k$ possible  
choices for their first index, we have in total at most   
$$
2^{C \log n_k} \, n_k^{C \log n_k} = (2 n_k)^{C \log n_k}
$$
possible choices. For each of these choices, it follows from (2.7) that at least
$n_k  (\theta \hv / (\theta - 1) - \eps')$ digits are prescribed (and are equal to $0$ or $b-1$). 
Consequently,  defining the positive real number $\eps''$ by the next equality, at most 
$$
n_k -n_k\Bigl({\theta \hv \over \theta - 1} - \eps' \Bigr) + 1
=n_k(1 + \eps'')\cdot  {\theta - 1 -  \theta \hv \over  \theta - 1} 
$$
digits in the prefix    
$a_1a_2 \ldots  a_{n_k}$, and thus in   
$a_1a_2 \ldots  a_{m_k}$, are free. 
The set of real numbers whose $b$-ary expansion starts with $a_1 a_2 \ldots  a_{m_k}$ 
defines an interval of length $b^{-m_k}$. 
By (2.5), we have 
$$
b^{-m_k}\leq b^{-(1 + \theta \hv) (1 - \eps'') n_k},
$$
for $k$ large enough. 
We have shown that the set of those $\xi$ corresponding to $(n_k)_{k \ge 1}, (m_k)_{k \ge 1}$ 
is covered by 
$$
(2n_k)^{C \log n_k} \,  b^{n_k(1 + \eps'') ( \theta - 1 -  \theta \hv ) / ( \theta - 1)} .
$$
intervals of length at most 
$ b^{-(1 + \theta \hv) (1 - \eps'') n_k}.$ 

Then, a standard covering argument shows that we have to consider the series
$$
\sum_{N \ge 1}  \, (2N)^{C \log N} \,  b^{N(1 + \eps'') ( \theta - 1 -  \theta \hv ) / ( \theta - 1)} 
\, b^{-(1 + \theta \hv) (1 - \eps'') N s}.    \eqno (2.9)
$$
The critical exponent $s_0$ such that (2.9) converges if $s > s_0$ 
and diverges if $s < s_0$ is given by
$$
s_0 = {1 + \eps'' \over 1 - \eps''} \cdot {\theta - 1 -  \theta \hv  \over (1 + \theta \hv) ( \theta - 1) }.
$$
It then follows that
$$
\dim (\{ \xi \in \R : \hv_b (\xi) \ge \hv\} \cap \{ \xi \in \R : v_b (\xi) = \theta \hv \}) \le 
{\theta - 1 -  \theta \hv  \over (1 + \theta \hv) ( \theta - 1) }.  \eqno (2.10)
$$
Actually, we have proved that, for every $\theta \ge  (1 - \hv)^{-1}$
and for every sufficiently small positive number $\delta$, we have
$$
\dim (\{ \xi \in \R : \hv_b (\xi) \ge \hv\} \cap 
\{ \xi \in \R :  \theta   \le v_b (\xi) \le  \theta + \delta  \}) \le 
{\theta - 1 -  \theta \hv  \over (1 + \theta \hv) ( \theta - 1) } + 5 \delta \hv.
$$

Regarding the right-hand side of (2.10) as a function of $\theta$ and taking (2.8)
into account, a short calculation shows that the maximum
is attained for $\theta = 2 / (1 - \hv)$, giving, for any positive $\eps$, that
$$
\dim \{ \xi \in \R : \hv_b (\xi) \ge \hv\}  
\le {1 + \eps \over 1 - \eps} \, {(1 - \hv)^2 \over (1 + \hv)^2}.
$$
We have established the required upper bound.

\vskip 4mm

\noindent{\bf 2.2. Lower bound.}

To obtain the lower bound, we construct a suitable Cantor type set. 
Let $\hv$ be in $(0,1)$. 
Let $\theta$ be a real number with $\theta\geq {1 \over 1-\hv}$. 
Choose two sequences $(m_k)_{k\geq 1}$ and $(n_k)_{k\geq 1}$ such that 
$n_k< m_k< n_{k+1}$ for $k \ge 1$, 
and such that $(m_k - n_k)_{k \ge 1}$ is non-decreasing. 
Furthermore, we assume that
$$
\lim_{k \to + \infty} \, {m_k - n_k \over n_{k+1}}=\hv, \eqno (2.11)
$$
and 
$$
 \lim_{k \to + \infty} \, {m_k - n_k \over n_k} = \theta \hv. \eqno (2.12)
$$
An easy way to construct such sequences is to start with
$$
n'_k= \left\lfloor \theta^k \right\rfloor, \quad m'_k
=\left\lfloor (\theta \hv+1) n'_k \right\rfloor, 
$$
and then to make a small adjustment to guarantee that $(m_k - n_k)_{k \ge 1}$ is 
non-decreasing.

We consider the set of real numbers $\xi$ in $(0, 1)$ whose $b$-ary expansion 
$\xi=\sum_{j \ge 1} {a _j \over b^j}$ satisfies, for $k \ge 1$, 
$$
 a_{n_k}=1, \ a_{n_k + 1} = \ldots = a_{m_k - 1} = 0, \ a_{m_k}=1, 
$$
and  
$$
\ a_{m_k+(m_k-n_k)}=  a_{m_k+2(m_k-n_k)}= \ \cdots \ = a_{m_k+t_k(m_k-n_k)}=1,
$$
where $t_k$ is the largest integer such that $m_k+t_k(m_k-n_k)< n_{k+1}$. Observe that, since
$$
t_k < {n_{k+1} - m_k \over m_k - n_k} \le {2 \over \hv},
$$
for $k$ large enough, the sequence $(t_k)_{k \ge 1}$ is bounded.

We check that the maximal length of blocks of zeros in the prefix of length $n_{k+1}$
of the infinite sequence $a_1 a_2 \ldots $ is equal to $m_k-n_k-1$. Thus, we deduce that
$$
\hv_b (\xi) =\hv \quad \hbox{and} \quad v_b(\xi) = \theta \hv.
$$
Actually, the above two equalities might not be true if $b=2$. 
However, they will be valid if we take the block $10$ in place of $1$ 
in the definition of $a_{m_k+t(m_k-n_k)}$ for $t = 1, \ldots , t_k$. 
The following proof will be almost the same. So, for simplicity, 
we assume that $b\geq 3$ and leave to the reader the slight change 
to deal with the case $b=2$ (or the reader may look at Subsection 4.2). 

Our Cantor type subset $E_{\theta, \hv}$ consists precisely of the real numbers in $(0, 1)$ 
whose $b$-ary expansion has the above property.
We will now estimate the Hausdorff dimension of $E_{\theta, \hv}$ from below.

%

Let $n$ be a large positive integer. 
For $a_1, \ldots , a_n$ in $\{0, 1, \ldots , b-1\}$,
denote by 
$
I_n(a_1, \dots, a_n )
$
the interval composed of the real numbers in $(0, 1)$ whose 
$b$-ary expansion starts with $a_1  \ldots a_n$.
Define a Bernoulli measure $\mu$ on $E_{\theta, \hv}$ as follows. 
We distribute the mass uniformly. 
For $k \ge 1$, set
$$ 
\delta_k:= m_k-n_k-1 \quad {\rm and} \quad u_k:=m_k+t_k(m_k-n_k).
$$

If there exists $k \ge 2$ such that $n_{k}\leq n\leq m_{k} $, then define
$$
\mu(I_n(a_1, \dots, a_n ))=b^{-(n_1-1+\sum_{j=1}^{k-1}(t_j\delta_j +n_{j+1}-u_j-1))}
=b^{-(n_k-1-\sum_{j=1}^{k-1}(m_j-n_j+1+t_j))}. 
$$ 
Observe that we have $\mu(I_{n_k})= \mu(I_{n_{k}+1}) = \ldots = \mu(I_{m_k})$.   

If there exists $k \ge 2$ such that $m_k< n<n_{k+1}$, then define 
$$
\eqalign{
\mu(I_n(a_1, \dots, a_n ))
&=b^{-(n_1-1+\sum_{j=1}^{k-1}(t_j\delta_j +n_{j+1}-u_j-1)+t\delta_k+n-(m_k+t\delta_k))}\cr 
&=b^{-n+\sum_{j=1}^{k-1}(m_j-n_j+1+t_j)+m_k-n_k+1+t}, }
$$
where $t$ is the largest integer such that $m_k+t(m_k-n_k)\leq n$.
It is routine to check that $\mu$ is well defined on $E_{\theta, \hv}$.

Now we calculate the {\it local dimension} of $\mu$ at $x\in E_{\theta, \hv}$, i.e., 
$$
\liminf_{r\to 0} {\log\mu(B(x,r)) \over \log r },
$$
where $B(x,r)$ stands for the ball centered at $x$ with radius $r$.
To this end, we first calculate the same lower limit for the basic intervals, and we prove that
$$
\liminf_{n\to \infty} {\log\mu(I_{n}) \over \log |I_{n}| } 
= {\theta-1-\theta\hv \over (\theta-1)(\theta\hv+1)}. \eqno (2.13)
$$
Since the lengths of basic intervals decrease `regularly', the above limit (2.13)   
is the same as the local dimension of $\mu$ at $x\in E_{\theta, \hv}$. 
The details are left to the reader. Finally, the required lower bound 
follows from the mass distribution principle; see \cite{Fal97}, page 26. 

Let us finish the proof by showing (2.13). First, we check that if $n=m_k$, then 
$$
\eqalign{
\liminf_{k\to \infty} {\log\mu(I_{m_k}) \over \log |I_{m_k}| }
&= \liminf_{k\to\infty} {n_k-1-\sum_{j=1}^{k-1}(m_j-n_j+1+t_j) \over m_k} \cr
&=\liminf_{k\to\infty}{n_{1}-1+\sum_{j=1}^{k-1}(n_{j+1}-m_j+1+t_j) \over m_{k}}.
}
$$
Recalling that $(t_k)_{k \ge 1}$ is bounded and that $(m_k)_{k \ge 1}$
grows exponentially fast in terms of $k$, we have 
$$
\liminf_{k\to \infty} {\log\mu(I_{m_k}) \over \log |I_{m_k}| }=
\liminf_{k\to\infty}{\sum_{j=1}^{k-1}(n_{j+1}-m_j) \over m_{k}}.$$
By (2.11) and (2.12), we see that
$$
\lim_{k\to\infty} {m_k \over n_k}= \theta\hv +1, \quad \lim_{k\to\infty} {m_{k+1} \over m_k}
= \theta, \quad {\rm and} \quad \lim_{k\to\infty} {n_{k+1} \over m_{k}}
=  {\theta \over \theta\hv +1}.
$$
Thus,
by the Stolz--Ces\`aro Theorem, 
$$
\eqalign{
\lim_{k\to\infty}{\sum_{j=1}^{k-1}(n_{j+1}-m_j) \over m_{k}} 
& = \lim_{k\to\infty}{n_{k+1}-m_{k}  \over m_{k+1}- m_{k}} \cr
& =  \lim_{k\to\infty}{ {n_{k+1} \over m_k}-1  \over {m_{k+1} \over m_k} - 1}
={\theta-1-\theta\hv \over (\theta-1)(\theta\hv+1)}. \cr}
$$
Hence, 
$$
\liminf_{k\to \infty} {\log\mu(I_{m_k}) \over \log |I_{m_k}| } 
={\theta-1-\theta\hv \over (\theta-1)(\theta\hv+1)}.
$$

Let $n$ be a large positive integer. 
If there exists $k \ge 2$ such that $n_{k}\leq n\leq m_{k} $, then 
$$
{\log\mu(I_{n}) \over \log |I_{n}| } \geq {\log\mu(I_{n}) \over \log |I_{m_k}| } 
= {\log\mu(I_{m_k}) \over \log |I_{m_k}| }.
$$
If there exists $k \ge 2$ such that $m_{k}<n<n_{k+1}$, 
write $n=m_k+t(m_k-n_k)+\ell$, where $t, \ell$
are integers with $0\leq t\leq t_k$ and $0\leq \ell<m_k-n_k$. Then we have 
$$
\mu(I_{n})=\mu(I_{m_k})\cdot b^{-(t\delta_k+\ell)} \quad {\rm and}
\quad |I_n|=|I_{m_k}|\cdot b^{-(t(m_k-n_k)+\ell)}.
$$
Since $0 \le t \le t_k$ and $(t_k)_{k \ge 1}$ is bounded, for $n$ large enough,
$$
{-\log\mu(I_{n}) \over -\log |I_{n}| }
= {-\log\mu(I_{m_k})+ t\delta_k+\ell \over -\log |I_{m_k}| + t(m_k-n_k)+\ell} 
\geq  {-\log\mu(I_{m_k}) \over -\log |I_{m_k}| },
$$
where we have used the fact that 
$$
{a+x \over b+x} \geq {a \over b}, \quad \hbox{for all} \ 0<a\leq b, x\geq 0.
$$
We have established (2.13).



\vskip 5mm

\centerline{\bf 3. Classical results on $\beta$-expansions}

\vskip 6mm
%

Throughout this section, $\beta$ denotes a real number greater than $1$
and $\lfloor \beta \rceil$ is equal to $\beta - 1$ if $\beta$ is an integer and
to $\lfloor \beta \rfloor$ otherwise.
The notion of $\beta$-expansion was introduced by
R\'enyi \cite{Re57} in 1957. 
We denote by $T_{\beta}$ the transformation defined
on $[0, 1]$ by $T_{\beta} (x) = \{\beta x\}$,
where $\{ \cdot \}$ denotes the fractional part function.  

\proclaim Definition 3.1. 
The expansion of a number $x$ in $[0, 1]$ to base $\beta$,
also called the $\beta$-expansion of $x$, 
is the sequence $(\eps_n)_{n \ge 1}=(\eps_n(x,\beta))_{n \ge 1}$ of integers 
from $\{0, 1, \ldots , \lfloor \beta \rceil\}$ such that
$$
x = {\eps_1 \over \beta} + {\eps_2 \over \beta^2} + \ldots 
+ {\eps_n \over \beta^n} + \ldots ,
$$
and, unless $x=1$ and $\beta$ is an integer,    
defined by one of the following equivalent properties:
$$
\sum_{k > n} \, {\eps_k \over \beta^k} < {1 \over \beta^n},
\quad \hbox{for all $n \ge 0$};  
$$
$$
\eps_1 = \lfloor \beta x \rfloor, \ \eps_2 = \lfloor \beta\{\beta x\} \rfloor, 
\ \eps_3 = \lfloor \beta \{\beta \{\beta x\} \} \rfloor, \ldots  
$$
$$
\eps_n = \lfloor \beta T_{\beta}^{n-1} (x) \rfloor, 
\quad \hbox{for all $n \ge 1$}.  
$$
We then write
$$
d_{\beta} (x) =  \eps_1 \eps_2 \ldots \eps_n \ldots.
$$

For $x<1$, the $\beta$-expansion coincides with the representation
of $x$ computed by the `greedy algorithm'.
If $\beta$ is an integer, then the digits $\eps_i$ of $x$ 
lie in the set $\{0,1,\ldots, \beta-1\}$ and $d_{\beta}(x)$ 
corresponds, for $x \not= 1$, 
to the usual $\beta$-ary expansion of $x$. If the $\beta$-expansion 
$$
d_\beta(1)=\eps_1(1, \beta)\eps_2(1, \beta)\dots \eps_n(1, \beta)\dots
$$ 
of $1$ is finite, i.e., if there exists $m\geq 1$ 
such that $\eps_m(1, \beta)\neq 0$ and $\eps_n(1, \beta)=0$ for all $n>m$, 
then $\beta$ is called a {\it simple Parry number}. In this case, 
we define the {\it infinite $\beta$-expansion of $1$} by 
$$
(\eps_1^*(\beta),\eps_2^*(\beta),\dots, \eps_n^*(\beta),\dots)
:=(\eps_1^*(1, \beta),\eps_2^*(1, \beta),\dots, (\eps_m(1, \beta)-1))^\infty,
$$
where $(w)^\infty$ stands for the periodic sequence $(w,w,w, \dots)$.

We endow the set 
$\{0, 1, \ldots , \lfloor \beta \rceil \}^{\Z_{\ge 1}}$
with  the product topology
and the one-sided shift operator $\sigma$ defined
by $\sigma ((s_n)_{n \ge 1}) = (s_{n+1})_{n \ge 1}$, for
any infinite sequence $(s_n)_{n \ge 1}$ in 
$\{0, 1, \ldots , \lfloor \beta \rceil \}^{\Z_{\ge 1}}$.

The lexicographic order on $\{0, 1, \ldots , \lfloor \beta \rceil \}^{\Z_{\ge 1}}$,  
denoted by $<_{\rm lex}$, is defined as follows: we write  
$$
w=(w_1,w_2,\dots ) <_{\rm lex} w'=(w'_1,w'_2,\dots ) 
$$
if there exists $k\geq 1$ such that for all $j<k$ we have $w_j=w'_j$, 
but $w_k<w'_k$. We use the notation $w\leq_{\rm lex}w'$ if $w<_{\rm lex}w'$ or $w=w'$.

\proclaim Definition 3.2.
The closure of the set of all $\beta$-expansions
of $x$ in $[0, 1]$ is called the $\beta$-shift  
and denoted by $\Sigma_{\beta}$.

Parry \cite{Pa60} proved that the $\beta$-shift $\Sigma_{\beta}$
is fully determined by $d_{\beta} (1)$.

\proclaim Theorem 3.3.
If $d_{\beta} (1) = \eps_1 \ldots \eps_m 00 \ldots 0 \ldots$ is finite with $\eps_m\neq 0$, then 
${\bf s} = (s_n)_{n \ge 1}$ belongs to $\Sigma_{\beta}$
if, and only if,
$$
\sigma^k ({\bf s}) \leq_{\rm lex} (\eps_1, \ldots, \eps_{m-1}, (\eps_m - 1) )^\infty, \quad
\hbox{for $k \ge 1$}.
$$
If $d_{\beta} (1)$ does not terminate with zeros only, then 
${\bf s} = (s_n)_{n \ge 1}$ belongs to $\Sigma_{\beta}$
if, and only if,
$$
\sigma^k ({\bf s}) \leq_{\rm lex} d_{\beta} (1), \quad
\hbox{for $k \ge 1$}.
$$

It follows from Theorem 3.3 that $\Sigma_{\beta}$ is contained in 
$\Sigma_{\beta'}$ if, and only if, $\beta \le \beta'$.

\proclaim Definition 3.4.
A block $\eps_1 \ldots \eps_m$, respectively, an infinite sequence   
$\eps_1 \eps_2 \ldots $, on the alphabet $\{0, 1, \ldots , \lfloor \beta \rceil\}$ 
is $\beta$-admissible (or, simply, admissible) if
$$
\sigma^k (\eps_1 \ldots \eps_m) \leq_{\rm lex} d_{\beta} (1), \quad
\hbox{for $k = 0, 1, \ldots , m-1$},
$$
respectively, if
$$
\sigma^k (\eps_1 \eps_2 \ldots ) \leq_{\rm lex} d_{\beta} (1), \quad
\hbox{for $k \geq 0$}.
$$
An infinite sequence $(\eps_1, \eps_2, \ldots)$ is self-admissible if
$$
\sigma^k (\eps_1 \eps_2 \ldots ) \leq_{\rm lex} (\eps_1 \eps_2 \ldots ), \quad
\hbox{for $k \ge 0$}.
$$

Denote by $\Sigma_\beta^n$ the set of all $\beta$-admissible blocks of length $n$. 
Then, its cardinality satisfies (R\'enyi \cite{Re57}, formula 4.9)   
$$
\beta^n \leq \sharp \Sigma_\beta^n \leq {\beta^{n+1} \over \beta-1}. \eqno (3.1)
$$
For any $(\eps_1, \dots, \eps_n)\in \Sigma_\beta^n$ , call 
$$
I_n(\eps_1, \dots, \eps_n):=\{x\in [0,1] : d_\beta(x) \ {\rm starts \ with} \  \eps_1 \ldots\eps_n\}
$$
an {\it $n$-th order basic interval} (with respect to the base $\beta$). 
Denote by $I_n(x)$ the $n$-th order basic interval containing $x$.
We remark that the basic intervals are also called cylinders by some authors. 

The next theorem was proved by Parry \cite{Pa60}.

\proclaim Theorem 3.5.
A sequence of digits $(\eps_1, \eps_2, \ldots)$ is the $\beta$-expansion of $1$
for some $\beta > 1$ if and only if it is self-admissible.

Now, we estimate the length of the basic intervals. 
We will use the notion of ``full cylinder" introduced
by Fan and Wang \cite{FaWa12}. 

\proclaim Definition 3.6. 
Let $(\eps_1, \dots, \eps_n)\in \Sigma_\beta^n$ 
be a $\beta$-admissible block in $\Sigma_\beta$. 
A basic interval $I_n(\eps_1, \dots, \eps_n)$  
is called full if it is of length $\beta^{-n}$.

\proclaim Proposition 3.7. (\cite{FaWa12}, Lemma 3.1)
A basic interval $I_n(\eps_1, \dots, \eps_n)$ is full 
if and only if, for any admissible block $(\eps'_1, \dots, \eps'_m)\in \Sigma_\beta^m$, 
the concatenation $(\eps_1, \dots, \eps_n, \eps'_1, \dots, \eps'_m)$ 
is also admissible, i.e., in $\Sigma_\beta^{n+m}$.

\proclaim Proposition 3.8. (\cite{ShWa13}, Corollary 2.6)
For any $w\in \Sigma_\beta^n$, if $I_n(w)$ is full, then for any $w'\in \Sigma_\beta^m$, we
have 
$$
|I_{n+m}(w, w')|=|I_n(w)|\cdot |I_m(w')|=\beta^{-n}|I_m(w')|.
$$

The following approximation of $\beta$-shift is very useful.
For $\beta>1$, let $(\eps_k^*(\beta))_{k\geq 1}$ denote the infinite  
$\beta$-expansion of $1$.
 Let $\beta_N$ be the unique real number which satisfies the equation 
 $$
 1={\eps_1^*(\beta) \over z} + \cdots +{\eps_N^*(\beta) \over z^N}. 
 $$
 Then we have $\beta_N<\beta$ and the sequence $(\beta_N)_{N \ge 1}$
 increases and converges to $\beta$ when $N$ tends to infinity.
 Furthermore, the subshift of finite type $\Sigma_{\beta_N}$ 
 is a subset of the $\beta$-shift $\Sigma_\beta$.  
 The subsets $\Sigma_{\beta_N}$ are increasing and converge to $\Sigma_\beta$.

\proclaim Proposition 3.9. (\cite{ShWa13}, Lemma 2.7)
For any $w\in \Sigma_{\beta_N}^n$ viewed as an element of $\Sigma_\beta^n$, we have
$$
\beta^{-(n+N)}\leq |I_{n}(w)|\leq \beta^{-n}.
$$

\vskip 5mm

\centerline{\bf 4. Proof of Theorem 1.4}

\vskip 6mm

\noindent {\bf 4.1. Upper bound.} 

The proof is essentially similar to that in Subsection 2.1.
Let $x \in [0, 1]$ be a real number and let
$$
x = {a_1 \over \beta} + {a_2 \over \beta^2} + \ldots 
$$
denote its $\beta$-expansion. Assume that $v_\beta (x)$ is positive.   

Define the increasing sequences $(n'_k)_{k \ge 1}$ and $(m'_k)_{k \ge 1}$
as follows: for $k \ge 1$, we have
$$
a_{n'_k} >0, a_{n'_k + 1} = \ldots = a_{m'_k - 1} = 0, a_{m'_k} > 0,
$$
and, furthermore, for every $j$ such that $a_j = 0$, there exists an index $k$
satisfying $n'_k < j < m'_k$.
Then take the maximal subsequences $(n_k)_{k \ge 1}$ and $(m_k)_{k \ge 1}$
of $(n'_k)_{k \ge 1}$ and $(m'_k)_{k \ge 1}$, respectively, in such a way that the sequence 
$(m_k - n_k)_{k \ge 1}$ is non-decreasing. Observe that, since $v_\beta (x) > 0$, this
sequence tends to infinity with $k$. 
Similarly, notice that
$$
\beta^{n_k - m_k} < T_\beta^{n_k} x < \beta^{n_k - m_k + 1}.
$$
We also have
$$
v_\beta (x)= \limsup_{k \to + \infty} \, {m_k - n_k \over n_k}, \quad
\hv_\beta (x) \le \liminf_{k \to + \infty} \, {m_k - n_k \over n_{k+1}},
$$
and the relations 
$$
\hv_\beta (x) \le  {v_\beta (x) \over 1 + v_\beta (x)} \quad {\rm{and}} \quad
v_\beta (x) \ge {\hv_\beta (x) \over 1 - \hv_\beta (x)}.
$$
The last inequality means that $v_\beta (x)$ is infinite when $\hv_\beta (x)= 1$. 
Combined with Theorem SW, it implies the last assertion of Theorem 1.4. 

Let $\hv$ be in $(0, 1)$. 
Let $\theta$ be the real number defined by
$$
v_\beta (x) = \theta \hv.
$$
Then
$$
\theta \ge {1 \over 1 - \hv}. 
$$
Arguing as in the proof of Theorem 1.2, we see that there
exists a positive real number $\eps$ such that the sum of all the lengths of the blocks of 
zeros in the prefix of length $n_k$ of the infinite sequence $a_1 a_2 \ldots $ is at least equal to
$$
n_k \, \Bigl( {\theta \hv \over \theta - 1} - \eps \Bigr),
$$
for $k \ge 1$.


Now we study the possible choices of digits among $a_1, \ldots , a_{m_k}$  
as we did in Subsection 2.1.
Note that there are $k$ blocks of digits which are `free'. 
Denote their lengths by $\ell_1, \dots, \ell_k$. Then we have $\ell_j=n_j-m_{j-1}$ and 
$$
\sum_{j=1}^k \ell_j \leq  n_k -n_k\left({\theta \hat{v} \over \theta-1}-\varepsilon\right)
={ n_k(1 + \eps') ( \theta - 1 -  \theta \hat{v} )/(\theta-1)},
$$
with $\eps' >0$ being a small real number. 
By (3.1), there are at most 
$$
{\beta \over \beta-1} \cdot \beta^{\ell_j} 
$$
ways to choose the block of length $\ell_j$. 
Thus, for the $k$ blocks, we have in total
$$
\left({\beta \over \beta-1}\right)^k \cdot \beta^{\sum_{j=1}^k\ell_j} 
\leq \left({\beta \over \beta-1}\right)^k \cdot \beta^{ n_k(1 + \eps') 
( \theta - 1 -  \theta \hat{v} )/(\theta-1)}
$$
choices.

As in Subsection 2.1, there exists a real number $C > 1$ 
such that $k$ is less than $C \log n_k$, and there are at most $n_k$ possible   
choices for the first index of the $k$ blocks. Thus, we have at most 
$$
n_k^{C \log n_k}
$$
possible choices for the $k$ blocks. 
We get that, for $k$ sufficiently large, the set of real numbers $x$   
with the above properties is contained in a union of no more than   
$$
\left({\beta n_k \over \beta-1}\right)^{C \log n_k} \beta^{ n_k(1 + \eps')
 ( \theta - 1 -  \theta \hv )/(\theta-1)}
$$
basic intervals of order $m_k$, each of them being of length at most $\beta^{-m_k}$.   

Furthermore, there are infinitely many indices $k$ such that
$$
\beta^{-m_k}\leq \beta^{-(1 + \theta \hv) (1 - \eps') n_k}.
$$
Arguing as in Subsection 2.1, we end up with the series
$$
\sum_{N \ge 1} \, N^{C \log N} \,  \beta^{N(1 + \eps') ( \theta - 1 -  \theta \hv )  / ( \theta - 1)} 
\, \beta^{-(1 + \theta \hv) (1 - \eps') N s}.  \eqno (4.1)
$$
The critical exponent $s_0$ such that (4.1) converges if $s > s_0$ 
and diverges if $s < s_0$ is given by
$$
s_0 = {1 + \eps'  \over 1 - \eps'} \cdot {\theta - 1 -  \theta \hv  \over (1 + \theta \hv) ( \theta - 1) }.
$$
We then get that
$$
\dim (\{ x \in (0, 1) : \hv_\beta (x) \ge \hv\} \cap \{ x \in (0, 1) : v_\beta (x) = \theta \hv \}) \le 
{\theta - 1 -  \theta \hv  \over (1 + \theta \hv) ( \theta - 1) }.   
$$
For the rest of the proof, we argue exactly as at the end of Subsection 2.1. We omit the
details. 

\vskip 4mm

\noindent {\bf 4.2. Lower bound.} 

Let $\hv$ be in $(0,1)$. 
Let $\theta$ be a real number with $\theta\geq {1 \over 1-\hv}$. 
We construct a subset $E_{\theta, \hv}$ whose elements $x$ satisfy   
$$
\hv_\beta (x)=\hv \quad \hbox{and} \quad v_\beta (x) = \theta \hv.
$$
We will suitably modify the construction we performed in Subsection 2.2 
when dealing with $b$-ary expansions.  

Let $N$ be a positive integer. 
Let $\beta_N$ be the real number defined from the infinite $\beta$-expansion of $1$   
as explained at the end of Section 3.  
As in Subsection 2.2, let $(m_k)_{k\geq 1}$ and $(n_k)_{k\geq 1}$ 
be sequences of positive integers with 
$n_k< m_k< n_{k+1}$ for $k \ge 1$, 
such that $(m_k - n_k)_{k \ge 1}$ is non-decreasing and (2.11) and (2.12) are satisfied.
Start with the construction performed in Subsection 2.2
and replace the digit $1$ for $a_{n_k}$,  
$a_{m_k}$ and $a_{m_k+j(m_k-n_k)},  1 \le j \leq t_k,$ by the block $0^N10^N$.  
Fill other places by blocks in $\Sigma_{\beta_N}$. Thus,  
we have completed the modifications and have constructed the subset $E_{\theta, \hv}$.   

Since $N$ is fixed and $(t_k)_{k \ge 1}$ 
is bounded, we check that for every $x$ in $E_{\theta, \hv} $,  
$$
\hv_\beta (x)= \lim_{k \to + \infty} \, {m_k - n_k-1+2N \over n_{k+1}+(4k+2)N+ \sum_{j=1}^{k} 2Nt_j}
=  \lim_{k \to + \infty} \, {m_k - n_k \over n_{k+1}}=\hv, 
$$
and
$$
v_\beta (x)= \lim_{k \to + \infty} \, {m_k - n_k-1+2N \over n_k+(4k-2)N+ \sum_{j=1}^{k-1} 2Nt_j}
=  \lim_{k \to + \infty} \, {m_k - n_k \over n_k}=\theta \hv.
$$
Moreover, the sequence $d_{\beta}(x)$ (in $\Sigma_\beta$) is, by construction,  
also in $\Sigma_{\beta_N}$.

Let $n$ be a large positive integer. 
Denote by 
$
I_n(a_1, \dots, a_n )
$
the interval composed of the real numbers in $(0, 1)$ whose 
$\beta$-expansion starts with $a_1 \dots  a_n$.
We will define a Bernoulli measure $\mu$ on $E_{\theta, \hv}$. 
We distribute the mass uniformly when we meet a block in $\Sigma_{\beta_N}$ 
and keep the mass when we go through the positions where 
the digits are determined by the construction. 
Precisely, we can write down the first levels as follows.

If $n<n_1$, define 
$$
\mu(I_n)={1 \over \sharp \Sigma_{\beta_N}^{n}}.
$$
If $n_1\leq n \leq m_1+4N$, then take
$$
\mu(I_n)={1 \over \sharp \Sigma_{\beta_N}^{n_1-1}}.
$$
If there exists an integer $t$ such that $0 \le t \le t_1$ and
$$
m_1+4N+ t(m_1-n_1)+2Nt< n 
\leq \min \{n_2+4N+2Nt_1, m_1+4N+(t+1)(m_1-n_1)+2Nt\},
$$ 
then set 
$$
\mu(I_n)={1 \over \sharp \Sigma_{\beta_N}^{n_1-1}} 
\cdot {1 \over (\sharp \Sigma_{\beta_N}^{m_1-n_1-1})^t} 
\cdot {1 \over \sharp \Sigma_{\beta_N}^{n-(m_1+4N+ t(m_1-n_1)+2Nt)}}.
$$
If there exists an integer $t$ such that $0 \le t \le t_1-1$ and
$$
m_1+4N+ (t+1)(m_1-n_1)+2Nt< n \leq m_1+4N+(t+1)(m_1-n_1)+2N(t+1),
$$ 
then set 
$$
\mu(I_n)={1 \over \sharp \Sigma_{\beta_N}^{n_1-1}} 
\cdot {1 \over (\sharp \Sigma_{\beta_N}^{m_1-n_1-1})^{t+1}} .
$$

More generally, for $k \ge 2$, set
$$
l_k:=n_{k}+(4k-4)N+ \sum_{j=1}^{k-1} 2Nt_j, \quad
h_k:=m_{k}+4kN+ \sum_{j=1}^{k-1} 2Nt_j.
$$
and
$$
\delta_k:=m_k-n_k-1, \quad u_k:=h_k+ t_k(m_k-n_k)+2Nt_k.
$$
If $l_k \leq n\leq h_k$,
define
$$
\mu(I_n)={1 \over \sharp \Sigma_{\beta_N}^{n_1-1}}
\cdot {1 \over \prod_{j=1}^{k-1}(\sharp \Sigma_{\beta_N}^{\delta_j})^{t_j}
\sharp \Sigma_{\beta_N}^{l_{j+1}-u_j-1}}
(=\mu(I_{l_k})=\mu(I_{h_k})).
$$
If there exists an integer $t$ such that $0 \le t \le t_k - 1$ and
$$ 
h_k+ t(m_k-n_k)+2Nt<n \leq \min\{l_{k+1}, h_k + (t+1)(m_k-n_k)+2Nt\},
$$ 
then define
$$
\mu(I_n)=\mu(I_{h_k})\cdot {1 \over (\sharp \Sigma_{\beta_N}^{\delta_k})^t}
 \cdot {1 \over \sharp \Sigma_{\beta_N}^{n-(h_k+ t(m_k-n_k)+2Nt)}}.
$$
If there exists an integer $t$ such that $0 \le t \le t_k - 1$ and
$$
h_k + (t+1)(m_k-n_k)+2Nt<n \leq h_k + (t+1)(m_k-n_k)+2N(t+1), 
$$
then define
$$
\mu(I_n)=\mu(I_{h_k})\cdot {1 \over (\sharp \Sigma_{\beta_N}^{\delta_k})^{t+1}}.
$$

By the construction and Proposition 3.7, we deduce that $I_{h_k}$ 
is full and thus has length  $\beta^{-h_k}$. 
Then, by recalling that $N$ is fixed and $(t_k)_{k \ge 1}$ is bounded, we 
deduce from (3.1) that
$$
\eqalign{
\liminf_{k\to \infty} {\log\mu(I_{h_k}) \over \log |I_{h_k}| }
&=\liminf_{k\to \infty} 
{\log \sharp \Sigma_{\beta_N}^{n_1-1} 
+ \sum_{j=1}^{k-1} \big({t_j} \log \sharp \Sigma_{\beta_N}^{\delta_j} 
+ \log \sharp \Sigma_{\beta_N}^{l_{j+1}-u_j-1} \big) \over h_k\log \beta } \cr
&= \lim_{k\to\infty}{n_1-1+\sum_{j=1}^{k-1} (t_j\delta_j+l_{j+1}-u_j-1) \over h_{k}} 
\cdot {\log \beta_N \over \log \beta}\cr
&=  \lim_{k\to\infty}{n_1-1+\sum_{j=1}^{k-1} (l_{j+1}-h_j-2Nt_j-1) \over h_{k}} 
\cdot {\log \beta_N \over \log \beta}\cr
&=\lim_{k\to\infty}{\sum_{j=1}^{k-1} (n_{j+1}-m_j )\over m_{k}}
\cdot {\log \beta_N \over \log \beta}\cr
&={\theta-1-\theta\hv \over (\theta-1)(\theta\hv+1)}\cdot {\log \beta_N \over \log \beta}.
}
$$

Let $n$ be a large positive integer.
If there exists $k \ge 2$ such that $l_k \leq n \leq h_k$, then 
 $$
{\log\mu(I_{n}) \over \log |I_{n}| } \geq {\log\mu(I_{n}) \over \log |I_{h_k}| } 
= {\log\mu(I_{h_k}) \over \log |I_{h_k}| }.
$$
If there exist integers $k \ge 2$ and $t$ such that $0 \le t \le t_k - 1$ and
$$
h_k+ t(m_k-n_k)+2Nt<n \leq \min\{l_{k+1}, h_k + (t+1)(m_k-n_k)+2Nt\},  
$$  
then, letting $\ell =n-(h_k+t(m_k-n_k)+2Nt)$, we have 
$$
\mu(I_{n})\leq \mu(I_{h_k})\cdot \beta_N^{-\delta_k t- \ell }.
$$

Since $I_{h_k}$ is full, then 
$$
 |I_n|=|I_{h_k}|\cdot |I_{n-h_k}(w')|,
$$
where $w'$ is an admissible block in $\Sigma_{\beta_N}$ of length $n-h_k$.
By Proposition 3.9,  
$$
 |I_{n-h_k}(w')|\geq \beta^{-(n-h_k+N)}.
$$
Hence,
$$
 |I_n|\geq |I_{h_k}|\cdot \beta^{-(n-h_k+N)}.
$$
Notice that 
$$
n-h_k+N=t(m_k-n_k)+2Nt+\ell+N=\delta_kt + \ell +t+ N(2t+1).
$$
Since $N$ is fixed and $t$ is bounded, we argue as in Subsection 2.2 to show that,    
for $n$ large enough, we have
$$
\eqalign{
{-\log\mu(I_{n}) \over -\log |I_{n}| } 
& \geq  {-\log\mu(I_{h_k})+(t\delta_k+\ell) \log \beta_N  \over 
-\log |I_{h_k}| +(t\delta_k+\ell +t+ N(2t+1))\log \beta} \cr
& \geq  {-\log\mu(I_{h_k}) \over -\log |I_{h_k}| }\cdot \eta(N), \cr}
$$
where $\eta(N)<1$ and $\eta(N)$ tends to $1$ as $N$ tends to infinity.   

Similarly, if there exist integers $k \ge 2$ and $t$ such that $0 \le t \le t_k - 1$ and
$$
h_k + (t+1)(m_k-n_k)+2Nt<n \leq h_k + (t+1)(m_k-n_k)+2N(t+1), 
$$
we also have for $n$ large enough,
$$
{-\log\mu(I_{n}) \over -\log |I_{n}| }\geq   {-\log\mu(I_{h_k}) \over -\log |I_{h_k}| }\cdot \eta(N).
$$
So, in all cases, we have
$$
\liminf_{n\to \infty} {\log\mu(I_{n}) \over \log |I_{n}| } 
\geq {\theta-1-\theta\hv \over (\theta-1)(\theta\hv+1)}\cdot {\log \beta_N \over \log \beta}\cdot \eta(N).
$$

Now, we consider a general ball $B(x,r)$ with $x$ a point    
in the Cantor-type subset   
$E_{\theta,\hv}$ and $r$ satisfying
$$
|I_{n+1}(x)|\leq r < |I_n(x)|.
$$
By the construction and Proposition 3.9, any $n$-th order basic interval $I_n$ satisfies     
$$
|I_n| \geq \beta^{-(n+N)}.
$$
Thus, the ball $B(x,r)$ intersects at most $\lceil 2 \beta^N\rceil +2 $ basic intervals of order $n$.
Noting that all $n$-th order basic intervals have the same measure, we deduce that  
$$
\mu(B(x,r)) \leq (\lceil 2 \beta^N\rceil +2)\cdot \mu(I_n(x)).
$$
On the other hand, by Proposition 3.9, 
$$
r\geq |I_{n+1}(x)| \geq 
\beta^{-(n+1+N)}=\beta^{-(N+1)}\cdot \beta^{-n} \geq \beta^{-(N+1)} \cdot |I_{n}(x)|
$$
Since $N$ is fixed, we have 
$$
\liminf_{r\to 0}{\log\mu(B(x,r)) \over \log r } 
=\liminf_{n\to \infty} {\log\mu(I_{n}(x)) \over \log |I_{n}(x)| }.
$$
Finally, the lower bound follows by letting $N$ tend to infinity. 

To end this section, we remark that the last step of the proof of the lower bound, i.e., 
the fact that the two lower limits concerning the basic intervals 
and the general balls are the same, always holds in the setting of $\beta$-transformation. 
This conclusion was proved very recently by Bugeaud and Wang in \cite{BuWa14}, who called
it a {\it modified mass distribution principle}. 

\vskip 5mm

\centerline{\bf 5. Proof of Theorem 1.6}

\vskip 6mm

We follow the approach of Persson and Schmeling \cite{PeSc08}. 
The main idea is to take a correspondence between the $\beta$-shift and the parameter space.
Then the results in the shift space can be translated to the parameter space.

\vskip 4mm

\noindent {\bf 5.1. Upper bound.}
  
First we reduce the question to a small interval $(\beta_0, \beta_1)$,
where $1 < \beta_0 < \beta_1$. 
For $\hv \in [0, 1]$ and $\theta \ge 1$, set     
$$
D_{\theta, \hv}:=\{ \beta > 1  : \hv_\beta (1) = \hv\} \cap \{ \beta>1: v_\beta (1) = \theta \hv \},
$$
and 
$$
D_{\theta, \hv}(\beta_0, \beta_1):=\{ \beta\in (\beta_0, \beta_1) : \hv_\beta (1) = \hv\} 
\cap \{ \beta \in (\beta_0, \beta_1): v_\beta (1) = \theta \hv \}.
$$
By Theorem 3.5, to every self-admissible sequence corresponds a real number $\beta > 1$. 
Let $K_{\beta_1}$ be the set of all self-admissible sequences in $\Sigma_{\beta_1}$. 
Let $\pi_{\beta}$ be the natural projection from the $\beta$-shift to the unit interval $[0,1]$. 
Then there is a one-to-one map 
$\varrho_{\beta_1} : \pi_{\beta_1}(K_{\beta_1}) \rightarrow (1, \beta_1)$.

Let $B_{\theta, \hv}$ be the subset of $\Sigma_{\beta_1}$ defined by   
$$
B_{\theta,\hv}:=\pi_{\beta_1}^{-1}\big(\{ x \in [0,1] : \hv_{\beta_1} (x) = \hv\} 
\cap \{ x \in [0,1] : v_{\beta_1} (x) = \theta \hv \}\big). 
$$
The H\"older exponent of the restriction of the map $\varrho_{\beta_1}$ to the subset 
$\pi_{\beta_1}(K_{\beta_1} \cap B_{\theta, \hv})$ is equal to  $\log \beta_0 / \log \beta_1$.
Since
$D_{\theta, \hv}(\beta_0, \beta_1)
\subset \varrho_{\beta_1}(\pi_{\beta_1}(K_{\beta_1} \cap B_{\theta, \hv}))$
we have
$$
\dim D_{\theta, \hv}(\beta_0, \beta_1) 
\leq \dim \varrho_{\beta_1}(\pi_{\beta_1}(K_{\beta_1} 
\cap B_{\theta, \hv})) \leq {\log \beta_1 \over \log \beta_0} \dim \pi_{\beta_1}(B_{\theta, \hv}),
$$
while, by Theorem 1.5, 
$$
 \dim \pi_{\beta_1}(B_{\theta, \hv}) =
{\theta - 1 -  \theta \hv  \over (1 + \theta \hv) ( \theta - 1) }.    
 $$
Letting $\beta_1$ tend to $\beta_0$, we obtain the requested upper bound.

\vskip 4mm

\noindent {\bf 5.2. Lower bound.} 

Take $\beta_2$ such that $1 < \beta_0<\beta_1< \beta_2$ 
and that the $\beta_2$-expansion of $1$ ends with zeros, i.e., such that the $\beta$-shift
$\Sigma_{\beta_2}$ is a subshift of finite type. We establish the following lemma.

\proclaim Lemma 5.1.
For real numbers $\beta_0, \beta_1, \beta_2, \hv$ and $\theta$ such that   
$1 < \beta_0<\beta_1< \beta_2$, $\hv \in [0, 1]$ and $\theta \ge 1$, we have 
$$
\dim(\varrho_{\beta_2}^{-1}(D_{\theta, \hv}(\beta_0, \beta_1)))
\geq {\theta - 1 -  \theta \hv  \over (1 + \theta \hv) ( \theta - 1) } 
\cdot {\log \beta_1 \over \log \beta_2}.
$$

It follows from Lemma 5.1 and the proof of Theorem 14 of 
Persson and Schmeling \cite{PeSc08}
that
$$
\dim(D_{\theta, \hv}(\beta_0, \beta_1))
\geq {\theta - 1 -  \theta \hv  \over 
(1 + \theta \hv) ( \theta - 1) } \cdot {\log \beta_1 \over \log \beta_2}.
$$
Finally, letting $\beta_2$ tend to $\beta_1$, we complete the proof of the lower bound.

\bigskip

\noindent {\it Proof of  Lemma 5.1.}

Denote by $(\varepsilon_k^*)_{k\geq 1}:=(\varepsilon_k^*(\beta_1))_{k\geq 1}$ 
the infinite $\beta_1$-expansion of $1$.  
Take an integer $N$ sufficiently large such that $\varepsilon_N^*\neq 0$ 
and the $\beta_0$-expansion of $1$ is smaller than $\eps_1^*\dots\eps_N^*$ 
in lexicographical order. 
Let $\tilde{\beta}_N$ be the unique real number which satisfies the equation 
$$
1={\eps_1^* \over z} + \cdots +{\eps_N^* \over z^N}. 
$$
Then we have $\beta_0<\tilde{\beta}_N<\beta_1$, and $\tilde{\beta}_N$ tends 
to $\beta_1$ as $N$ tends to infinity. 
Moreover, the infinite $\tilde{\beta}_N$-expansion of $1$ is given by
$$
(\eps_1^*,\dots, (\eps_N^*-1))^{\infty}.
$$
 
We construct a subset $\pi_{\beta_2}(K_N)$ of $\varrho_{\beta_2}^{-1}(\DD)$ 
by using the same construction as in Subsection 4.2. 
Take any sequence $\underline{a} \in \Sigma_{\tilde{\beta}_N}$ 
constructed in Subsection 4.2 with $\beta= \tilde{\beta}_N$.    
Make the concatenation $\eps_1^*\dots \eps_N^*0^N\underline{a}$.  
By Lemma 5.2 from \cite{LPWW}, the sequence $\eps_1^*\dots \eps_N^*0^N\underline{a}$ 
is self-admissible and thus, by Theorem 3.5,  
is the $\beta$-expansion of $1$ for some $\beta$. 
By checking the lexicographic ordering, 
we see that  $\beta_0 < \beta < \beta_1$.  
We define the subset $K_N$ 
to be the collection of these sequences. 
Notice that $K_N$ is also a subset of $\Sigma_{\beta_2}$. 

Now we can define a measure $\tilde{\mu}$ on the set $\pi_{\beta_2}(K_N)$.
Consider a basic interval $I_{2N+m}(\eps_1^*, \ldots , \eps_N^*, 0^N, a_1, \ldots , a_m)$. 
Let $\mu$ be the measure defined in Section 4.2 by replacing $\beta_N$ 
there by $\tilde{\beta}_N$. Define
$$
\tilde{\mu}(I_{2N+m}(\eps_1^*,\ldots, \eps_N^*,0^N,a_1, \ldots, a_m))
:=\mu(I_m(a_1, \ldots, a_m)).
$$ 

By Proposition 3.7, the basic interval $I_{2N}(\eps_1^*,\ldots, \eps_N^*,0^N)$ is full. 
Then, it follows from Proposition 3.8 that
$$
\eqalign{
|I_{2N+m}(\eps_1^*,\dots, \eps_N^*,0^N,a_1, \ldots, a_m)|
&=|I_{2N}(\eps_1^*,\dots, \eps_N^*,0^N)|\cdot |I_m(a_1, \ldots, a_m)|,\cr
&={\beta_2}^{-2N} |I_m(a_1, \ldots, a_m)|. \cr}
$$  
Since $\Sigma_{\beta_2}$ is of finite type, there exists a positive real number $C$ such that
$$
C^{-1} \beta_2^{-m} \leq |I_m(a_1, \ldots, a_m)|\leq C \beta_2^{-m}.
$$
Thus, noting that $N$ is fixed, we can deduce, as in Section 4.2, that 
$$
\liminf_{n\to \infty} {\log\mu(I_{n}) \over \log |I_{n}| } 
\geq {\theta-1-\theta\hv \over (\theta-1)(\theta\hv+1)}\cdot 
{\log \tilde{\beta}_N \over \log \beta_2}\cdot \eta(N).
$$
Similarly as in the previous section, we have the same inequality for the 
general ball $B(x,r)$, i.e., 
for any $x$ in the Cantor-type set $\pi_{\beta_2}(K_N)$, we have    
$$
\liminf_{r\to 0} {\log\mu(B(x,r)) \over \log r } 
\geq {\theta-1-\theta\hv \over (\theta-1)(\theta\hv+1)}
\cdot {\log \tilde{\beta}_N \over \log \beta_2}\cdot \eta(N).
$$
Hence
$$
\dim (\varrho_{\beta_2}^{-1}(\DD)) \geq \dim(\pi_{\beta_2}(K_N)) 
\geq {\theta-1-\theta\hv \over (\theta-1)(\theta\hv+1)}
\cdot {\log \tilde{\beta}_N \over \log \beta_2}\cdot \eta(N).
$$
Letting $N$ tend to infinity, this proves the lemma.   \cqfd

\vskip 5mm

\goodbreak

\centerline{\bf 6. Diophantine approximation and Cantor sets}

\vskip 6mm

Throughout this section, $b$ denotes an integer at least equal to $3$
and $S$ is a subset of $\{0, 1, \ldots , b-1\}$ of cardinality at least two and containing at least
one of the digits $0$ and $b-1$. Let $K_{b, S}$ denote the set of real numbers $\xi$
in $[0, 1]$ which can be expressed as
$$
\xi = \sum_{ j \ge 1} \, {a_j \over b^j},
$$
with $a_j \in K_{b, S}$ for $j \ge 1$. Note that for $b=3$ and $S = \{0, 2\}$ the set
$K_{b, S}$ is the middle third Cantor set, which we simply denote by $K$.
Furthermore, recall that the Hausdorff dimension of $K_{b, S}$ is given by
$$
\dim K_{b, S} = {\log \sharp S \over \log b}.
$$

As a corollary of a more general result, 
Levesley, Salp and Velani \cite{LeSaVe07} established that,
for $v \ge 0$,
$$
\dim \{\xi \in K : v_3 (\xi) \ge  v \} = 
\dim \{\xi \in K : v_3 (\xi) =  v \} = {\log 2 \over \log 3} \cdot {1 \over 1 + v}.   \eqno (6.1)
$$
A suitable modification of the proofs of Theorems 1.1 and 1.2 allows us
to extend (6.1) as follows.

\proclaim Theorem 6.1.
Let $b \ge 3$ be an integer  
and $S$ a subset of $\{0, 1, \ldots , b-1\}$ 
of cardinality at least two and containing at least
one of the digits $0$ and $b-1$.
Let $\theta$ and $\hv$ be positive real numbers with $\hv < 1$.
If $\theta < 1/(1 - \hv)$, then the set
$$
\{ \xi \in K_{b, S} : \hv_b (\xi) \ge \hv\} \cap \{ \xi \in K_{b, S} : v_b (\xi) = \theta \hv \}
$$
is empty. Otherwise, we have
$$
\dim (\{ \xi \in K_{b, S} : \hv_b (\xi) = \hv\} \cap \{ \xi \in K_{b, S} : v_b (\xi) = \theta \hv \}) =
{\log \sharp S \over \log b} \cdot 
{\theta - 1 -  \theta \hv  \over (1 + \theta \hv) ( \theta - 1) }.     \eqno (6.2)
$$
Furthermore,
$$
\dim \{ \xi \in K_{b, S} : \hv_b (\xi) =  1\} = 0.
$$
and
$$
\dim \{\xi \in K_{b, S} : \hv_b (\xi) \ge  \hv \} =
\dim \{\xi \in K_{b, S} : \hv_b (\xi) = \hv \} =
{\log \sharp S \over \log b} \cdot \Bigl( {1-\hv \over 1+\hv} \Bigr)^2.   
$$

The assumption that $S$ contains at least
one of the digits $0$ and $b-1$ is necessary, since, otherwise,
we trivially have $v_b (\xi) = 0$ for every $\xi$ in $K_{b, S}$. 

For additional results on Diophantine approximation on Cantor sets, the reader
may consult Chapter 7 of \cite{BuLiv2} and the references quoted therein.

\medskip

\noindent {\it Proof of Theorem 6.1.}  
We follow step by step the proof of Theorems 1.1 and 1.2. To prove that 
the right hand side of (6.2) is an upper bound for the dimension, we use the same 
covering argument, but, instead of (2.9), we have to consider the series
$$
\sum_{N \ge 1}  \, (2N)^{C \log N} \,  
(\sharp S)^{N(1 + \eps') ( \theta - 1 -  \theta \hv )  / ( \theta - 1)}  
\, b^{-(1 + \theta \hv) (1 - \eps') N s}.
$$
As for the lower bound, we again consider a Bernoulli measure and we distribute the
mass uniformly among the elements of $S$. Also, if $1$ does not belong to $S$, we 
cannot take $a_{n_k}, a_{m_k}, \ldots$ equal to $1$ and we then choose them equal to a
non-zero element of $S$. We omit the details. 
\cqfd

\bigskip
\noindent{\bf Acknowledgements.} 
L. Liao was partially supported by the ANR, grant 12R03191A -MUTADIS, France.

\vskip 8mm

\goodbreak

\centerline{\bf References}

\vskip 5mm

\beginthebibliography{999}

\bibitem{AmBu10} 
M. Amou and Y. Bugeaud,
{\it Expansions in integer bases and exponents of
Diophantine approximation},
J. London Math. Soc. 81 (2010), 297--316.

\bibitem{BeVe06}
V. Beresnevich and S. Velani, 
{\it A mass transference principle and the Duffin--Schaeffer conjecture for Hausdorff measures},
Ann. of Math. 164  (2006), 971--992.

\bibitem{BoFr72}
I. Borosh and A. S. Fraenkel,
{\it A generalization of Jarn\'\i k's 
theorem on Diophantine approximations},
Indag. Math. 34 (1972), 193--201.

\bibitem{BuLiv2}
Y. Bugeaud,
Distribution modulo one and Diophantine approximation.
Cambridge Tracts in Mathematics 193, Cambridge, 2012.

\bibitem{BuLa05a}
Y. Bugeaud and M. Laurent,
{\it Exponents of Diophantine Approximation and
Sturmian Continued Fractions},
Ann. Inst. Fourier (Grenoble) 55 (2005), 773--804.

\bibitem{BuWa14}
Y. Bugeaud and B.-W. Wang,
{\it Distribution of full cylinders and the Diophantine properties 
of the orbits in $\beta$-expansions},
J. Fractal Geometry. To appear.

 \bibitem{Fal97} K. Falconer, 
 Techniques in fractal geometry. John Wiley \& Sons, Ltd., Chichester, 1997.

\bibitem{FaWa12}
A.-H. Fan and B.-W. Wang,
{\it On the lengths of basic intervals in beta expansions},
 Nonlinearity,  25 (2012), 1329--1343.

\bibitem{Kh26}
A. Ya. Khintchine,
{\it \"Uber eine Klasse linearer diophantischer Approximationen},
Rendiconti Circ. Mat. Palermo 50 (1926), 170--195.

\bibitem{LeSaVe07}
J. Levesley, C. Salp and S. L. Velani,
{\it On a problem of K. Mahler: Diophantine
approximation and Cantor sets},
Math. Ann. 338 (2007), 97--118.

\bibitem{LPWW}
B. Li, T. Persson, B.-W. Wang and J. Wu,
{\it Diophantine approximation of the orbit of $1$ 
in the dynamical system of beta expansions},
Math. Z.   276 (2014), 799--827.

\bibitem{Pa60}
W. Parry, 
{\it On the $\beta$-expansions of real numbers}, 
Acta Math. Acad. Sci. Hung. 11 (1960), 401--416.

\bibitem{PeSc08}
T. Persson and J. Schmeling,
{\it Dyadic Diophantine approximation and Katok's horseshoe approximation}, 
Acta Arith 132 (2008), 205--230.

\bibitem{Re57}
A. R\'enyi,
{\it Representations for real numbers and their ergodic properties},
 Acta Math. Acad. Sci. Hung. 8 (1957), 477--493.

\bibitem{ShWa13}
L. M. Shen and B. W. Wang, 
{\it Shrinking target problems for beta-dynamical system}, 
Sci. China Math. 56 (2013), 91--104.

\endthebibliography

\vskip1cm

\noindent Yann Bugeaud  \hfill{Lingmin Liao  }

\noindent Universit\'e de Strasbourg    \hfill{ Universit\'e Paris-Est Cr\'eteil}

\noindent IRMA    \hfill{LAMA}

\noindent 7, rue Ren\'e Descartes      \hfill{61, av G\'en\'eral de Gaulle    }

\noindent 67084 STRASBOURG  (FRANCE)   \hfill{94000 CR\'ETEIL  (FRANCE)}

\vskip2mm

\noindent {\tt bugeaud@math.unistra.fr}    \hfill{{\tt lingmin.liao@u-pec.fr}}

\bye

----------------------------------------

\bibitem{BuLiv}
Y. Bugeaud,
Approximation by algebraic numbers,
Cambridge Tracts in Mathematics 160, Cambridge, 2004.

\bibitem{Bu08}
Y. Bugeaud,
{\it Diophantine approximation and Cantor sets},
Math. Ann. 341 (2008), 677--684.

\bibitem{BuLa07}
Y. Bugeaud and M. Laurent,
{\it Exponents of Diophantine approximation}.
In: Diophantine Geometry Proceedings,
Scuola Normale Superiore Pisa, Ser. CRM, vol. 4, 2007, 101--121.

\vskip 5mm

\centerline{\bf 6. Inhomogeneous approximation}

\vskip 6mm

We may easily extend Theorems 1.1 and 1.2 to inhomogeneous
approximation, by defining, for a given real number $\xi_0$ in $[0, 1]$
and for $\xi, b$ as in Definition 1.1, the quantities $v_b(\xi, \xi_0)$
and $\hv_b(\xi, \xi_0)$ to be the supremum of  the real numbers 
$v$ for which, for arbitrarily large (resp., for every sufficiently large)
integer $N$, the equation
$$
|| b^n \xi - \xi_0|| < (b^N)^{-v}
$$
has a solution $n$ with $1 \le n \le N$. 

Theorem SW has been extended to inhomogeneous approximation
by Shen and Wang \cite{ShWa13} and, subsequently, 
by Bugeaud and Wang \cite{BuWa14}.

\proclaim Theorem BW. 
Let $\beta > 1$ be a real number and $v$ be a positive real number. Then,
for any fixed sequence $\uy = (y_n)_{n \ge 1}$, we have
$$
\dim \{x \in (0, 1) : v_{\beta} (x, \underline{y}) \ge v\} = {1 \over 1 + v}.
$$

Theorem PS has been recently considerably extended
by Li, Persson, Wang and Wu \cite{LPWW}. 

\proclaim Theorem LPWW.
Let $\beta_0, \beta_1$ and $v$ be real numbers with $1 < \beta_0 < \beta_1$
and $v > 0$. Then, for every $x_0$ in $[0, 1]$, we have
$$
\dim \{\beta \in (\beta_0 , \beta_1) : v_{\beta} (1, x_0) \ge v \} = {1 \over 1 + v}.
$$
More generally, for every Lipschitz continuous function
$\beta \mapsto y(\beta)$ on $(1, + \infty)$ taking its values in $[0, 1]$,
we have
$$
\dim \{\beta \in (\beta_0 , \beta_1) : v_{\beta} (1, y (\beta)) \ge v \} = {1 \over 1 + v}.
$$

Combining our arguments with those from \cite{BuWa14,LPWW}, we can establish the
following results.

\proclaim Theorem 6.1.
Let $\beta > 1$ be a real number and $\hv$ be a real number in $(0, 1)$.
Then, we have
$$
\dim \{x \in (0, 1) : \hv_{\beta} (x, \uy) \ge \hv \} = \Bigl( {1- \hv \over 1+ \hv} \Bigr)^2.
$$

\proclaim Theorem 6.2. 
Let $\beta_0, \beta_1$ and $v$ be real numbers with $1 < \beta_0 < \beta_1$
and $v > 0$. Then, for every Lipschitz continuous function
$\beta \mapsto y(\beta)$ on $(1, + \infty)$ taking its values in $[0, 1]$,
we have
$$
\dim \{\beta \in (\beta_0 , \beta_1) : \hv_{\beta} (1, y(\beta)) \ge \hv \} 
= \Bigl( {1- \hv \over 1+ \hv} \Bigr)^2.
$$


\proclaim Theorem 1.7. 
Let $\beta_0, \beta_1$ be real numbers with $1 < \beta_0 < \beta_1$.
Let $\hv$ be a real number in $[0, 1]$.
Then we have
$$
\dim \{\beta \in (\beta_0 , \beta_1) : \hv_{\beta} (1) \ge \hv \} 
= \Bigl( {1- \hv \over 1+ \hv} \Bigr)^2.
$$

\bye